\documentclass[12pt]{elsarticle}
\usepackage{graphicx}
\usepackage{amssymb}
\usepackage{amsmath}
\usepackage{subfigure}

\newcommand\Rb{\mathbf{R}}
\newcommand\Bb{\mathbf{B}}
\newcommand\xb{\mathbf{x}}
\newcommand\yb{\mathbf{y}}
\newcommand\zetahat{\boldsymbol{\hat{\zeta}}}
\newcommand\Zhat{\mathbf{\hat{Z}}}
\begin{document}

\title{A partially mesh-free scheme for representing anisotropic spatial variations along field lines}

\author{Ben F McMillan}
\address{Centre for Fusion, Space and Astrophysics, University of Warwick, Coventry}
\date{\today}

\begin{abstract}
  A common numerical task is to represent functions which are highly spatially anisotropic, and to solve
  differential equations related to these functions. One way such anisotropy arises is that
  information transfer along one spatial direction is much faster than in others. In this situation,
  the derivative of the function is small in the local direction of a vector
  field $\Bb$. In order to define a discrete representation, a set of surfaces $M_i$ indexed by an integer $i$
  are chosen such that mapping along the field $\Bb$ induces a one-to-one relation between the points on surface
  $M_i$ to those on $M_{i+1}$. For simple cases $M_i$ may be surfaces of constant coordinate value.
  On each surface $M_i$, a function description is constructed using basis functions defined on 
  a regular structured mesh. The definition
  of each basis function is extended from the surface $M$ along the lines of the field $\Bb$ by multiplying it by a smooth compact
  support function whose argument increases with distance along $\Bb$. Function values are evaluated by summing
  contributions associated with each surface $M_i$. This does not require any special connectivity of the meshes
  used in the neighbouring surfaces $M$, which substantially simplifies the meshing problem compared to attempting to find
  a space filling anisotropic mesh. We explore the numerical properties of the scheme, and show that it can be used to efficiently
  solve differential equations for certain anisotropic problems.
\end{abstract}

\maketitle

\section{Introduction}

The technique proposed here is motivated by plasma physics examples where particles travel much more
easily along magnetic field lines than in the perpendicular direction, so that in quantities like fluid moments elongated structures are formed, aligned with the
field lines. In particular, the technique is designed to solve problems in magnetic confinement fusion (MCF),
where the field lines wind around a central axis and may be closed, trace out surfaces, or fill ergodic regions. An additional difficulty
in MCF problems is that the anisotropic structures are strongly curved, because field lines are not straight (even in cylindrical coordinates)
over the length scale of the structures; the departure from straightness is often considerably larger than the wavelength of the structure
in the directions of rapid variation. This paper outlines a method for representing functions aligned along field lines which are not necessarily
aligned on nested surfaces or closed (such as plasmas with an X-point), and for solving equations relating these functions.

A variety of techniques to deal with representing these highly anisotropic functions exist.
The canonical technique is to define a 3D mesh to fill the space of interest, with the mesh strongly elongated along the field
line. Achieving a very good alignment of the mesh along the field lines is in general quite a difficult meshing
problem, and for this reason many MCF physics codes work only in the region where the field lines trace out a
nested set of topologically toroidal magnetic surfaces: these are KAM tori\cite{KAM} associated
with the field line Hamiltonian. In the tokamak core, for example, because of near-axisymmetry, nested surfaces usually exist and
regular grids can efficiently be generated, or angular coordinates may be employed in conjunction with a Fourier representation.
This is not the case for stellarator geometry or in the tokamak edge region. 

To avoid difficult meshing problems for the general case where the region of interest is not filled by nested surfaces,
it is desirable to relax the requirement of mesh connectivity. The
{\it Flux Coordinate Independent} (FCI) approach\cite{Hariri2013,Ottaviani2011,Stegmeir}, based on a finite difference method, defines
function values on nodes lying on a set of surfaces $M_i$ which are taken to be surfaces of constant coordinate $\zeta$.
A node $\xb$ on surface $M_i$ can be mapped along the field direction $\Bb$ to find image points, $\xb_\pm$ on surfaces $M_{i \pm 1}$.
Although these image points will not in general lie on nodes on the surfaces $M_{i\pm 1}$, the function
may be evaluated at these points by interpolation. Given the values of the function at points $\xb_\pm$, derivatives along the
field direction may then be determined.

Another way to relax the mesh connectivity constraint is via a finite volume technique, where the volumes are extrusions of a
polygonal grid cell on one surface to the next, and a polynomial representation is chosen in each volume element;
smoothness constraints are then approximately imposed using a discontinuous Galerkin approach. A hybrid method incorporating
finite differences along the field line and the discontinuous Galerkin method has also been investigated\cite{Held201629}.

A natural method for representing the anisotropic functions of interest is to change coordinates by defining
a grid on a surface, and extending this to a volume grid by defining     
an additional coordinate parameterising the distance along the mapping (this is known as the flux tube method\cite{fluxtubes} in MCF).
Locally, this allows for straightforward and efficient representation of the problem anisotropy. However, the coordinate scheme
becomes highly distorted for mappings with strong shear or compression. The mesh connectivity problem also resurfaces if
the originating surface is eventually mapped back onto itself, as at this point the representation on two
non-aligned meshes must be combined in some way.

We propose a partially mesh-free method which we call FCIFEM as it is a Finite Element Method translation
of the FCI approach. The method represents anisotropic functions using a compact-support
set of basis functions which are defined in a local set of coordinates aligned with the mapping. The definition
of the set of basis functions is used to define weak forms of differential equations, as in a standard Galerkin method. 
The philosophy is to design a method which is robust and simple to implement, and requires little manual user interaction,
because it avoids complex mesh generation tasks. The representation also provides a simple way to neatly tackle a series of
related problems with slightly different configurations, generated, for example, when the field generating the mapping
varies slowly with time.

\section{Definition of the finite dimensional representation}

For the sake of simplicity, we consider a 3D volume labelled by coordinates $R$, $Z$, $\zeta$, and
take the surfaces $M$ of interest to be surfaces of constant coordinate $\zeta$, of value $\zeta_i$
on surface $M_i$.

Consider a continuous function $\mathcal{Q}: \mathbb{R}^3 \times \mathbb{R} \rightarrow \mathbb{R}^3$, that we will refer to as the {\it mapping}, which takes a
point $\xb$ and a parameter $s$ and returns a point $\yb$. We will use this function to define projections of the 3D space, along curves locally aligned with
the direction of the anisotropy, onto each of the surfaces $M_i$; we require $\mathcal{Q}(\xb,s)_{\zeta} = s$, so that the projection associated with surface $M_i$
has parameter $s=\zeta_i$. We also require $Q(\xb,x_\zeta)=\xb$ so points on the surface
map to themselves.

\begin{figure}[htb]
\centering
\includegraphics[width=8cm]{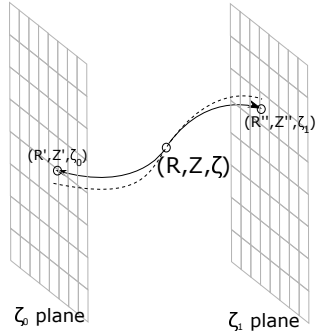}
\caption{ The geometry of the approximate flow line mapping. Here, the mapping from the point $\xb = (R,Z,\zeta)$ to the point
  $(R',Z',\zeta_0)$ on surface $M_0$ is depicted, with $(R',Z') = \mathcal{Q}(\xb,\zeta_0)$, as well as the analogous mapping
  to surface $M_1$. The dashed line shows the field line from the point $\xb$, for which eq. \ref{equation:mappingcond}
  holds exactly.
   }
\label{fig:mapping}
\end{figure}

One way to generate such a mapping would be to consider the action of a static flow field $\Bb$ displacing the position, leading to a field line equation
\begin{equation}
  \frac{d\xb}{dt} = \Bb(\xb).
\end{equation}
If we followed a field
line from position $\xb$ until it had toroidal coordinate $s$, where the field line was at the point $\yb$, we could define the mapping as
$\mathcal{Q}(\xb,s) = \yb$. We will refer to this as an exact mapping, which satisfies the equation
\begin{equation}
  \frac{\partial}{\partial \epsilon} \mathcal{Q}(\xb+\epsilon \Bb,\zeta)|_{\epsilon=0} = 0.
  \label{equation:mappingcond}
\end{equation}
It is convenient to allow the mapping function to be more general, however, and not necessarily exactly be the solution
to this field line mapping equation (or to any equation with a modified $\Bb$), either because we don't know the exact solution,
or because an approximate solution is numerically more
desirable. This has consequences for the quality of approximation, as the anisotropic direction will not align exactly with the mapping
direction, and certain statements on convergence will be shown only for the case of an exact mapping.
For consistency properties to hold, the mapping will be required to be one-to-one and at least of the same order of smoothness as the
element functions defined in the next paragraph. The geometry of this mapping is show in figure \ref{fig:mapping}.


\begin{figure}[htb]
\centering
\includegraphics[width=8cm]{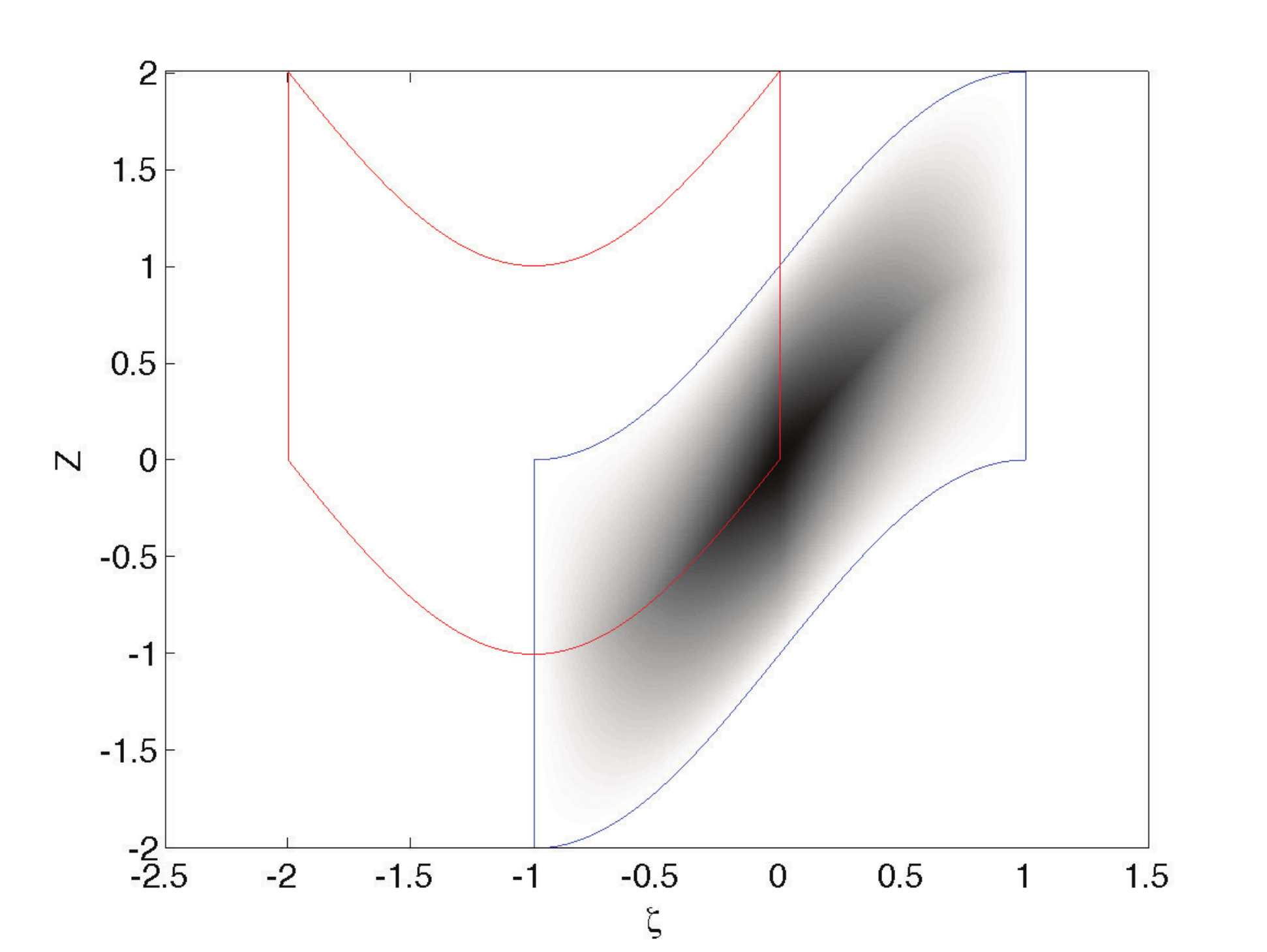}
\caption{ A density plot of the function $\phi$ for a 2D restriction of the FCIFEM representation to the plane $(\zeta,Z)$, with a single coefficient $\phi_{i,j,k}$ non-zero
  (equivalently, this can be seen as a slice at constant $R$ of a 3D FCIFEM representation). The domain of support of this nonzero term in the FCIFEM sum is shown as a blue
  line, and a domain of support of an element based on the plane $\zeta=-1$ is shown in red.
  The field $\Bb$ is such that field lines are of the form $ Z = \sin(\pi \zeta / 2) + Z_0$,
  and an exact mapping is used. Linear B-Spline basis functions are used for $\Omega_Z,\Omega_{\zeta}$.}
\label{fig:spline_function}
\end{figure}

The second element of the FCIFEM method to be chosen is the representation on the planes $M_i$. In general in might be helpful
to choose a general unstructured mesh, but
for the purposes of explanation and initial testing in this paper, we will use a simple uniformly spaced Cartesian mesh on each plane $M$.
In the interior region the representation of a scalar function of position is defined as
\begin{align}
  \phi(R,Z,\zeta) = \sum_{i,j,k} \phi_{i,j,k} &  \notag \\
  \times & \Omega_{R} \left[\mathcal{R}(R,Z,\zeta,\zeta_k) - R_i \right] \notag \\
  \times & \Omega_{Z} \left[\mathcal{Z}(R,Z,\zeta,\zeta_k) - Z_j \right]  \notag \\
  \times & \Omega_{\zeta} \left[ \zeta-\zeta_k \right]
  \label{eq:FCI_FEM}
\end{align}
with compact support basis functions $\Omega$, a regular set of Cartesian nodes $i,j,k$, and using the componentwise notation
$\mathcal{Q} = (\mathcal{R},\mathcal{Z},\zeta)$. We will choose the functions $\Omega$
to be B-Spline basis functions for the remainder as their properties are sufficient to
ensure smoothness and lowest order consistency (and this is similar to a finite element approach used
earlier in MCF codes\cite{Fivaz_1998,SebORB5}). An example of the shape of a distorted 2D basis function (the coefficient of $\phi_{i,j,k}$ for some chosen
$i,j$ and $k$) is plotted in figure \ref{fig:spline_function}. To evaluate the function value at point $\xb$, each term of the sum in eq. \ref{eq:FCI_FEM}
is evaluated by calculating the mapping $\mathcal{Q}(\xb,\zeta_k)$, and the product of the basis functions can then be directly calculated.
For a smooth mapping, the overall representation smoothness depends on the order of the spline. The arguments about convergence are most simply
made in the case with uniform nodes where $R_i = i \, \delta R $, $Z_j = j \, \delta Z $ and $\zeta_k = k \, \delta \zeta$. The space spanned
by these functions will be denoted $S$.

\begin{figure}[htb]
\centering
\includegraphics[width=8cm]{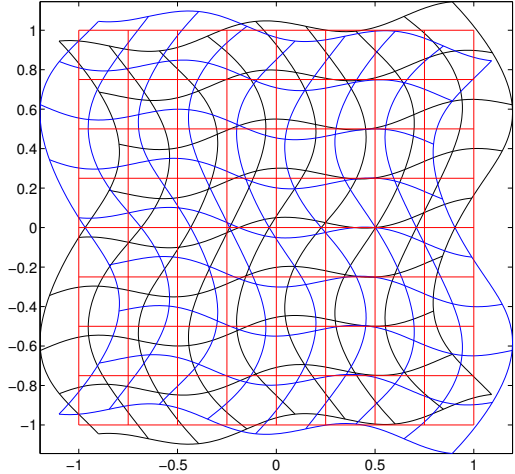}
\caption{ The intersection of mesh domains (which form the boundaries of domains of support of basis functions) with a constant $\zeta$ plane,
  for an example mapping. The red mesh is associated with a regular meshing on this plane $\zeta = \zeta_0$, so is rectilinear, whereas the blue and
  black meshes are associated with the planes $\zeta_{\pm 1}$. In general the domains of intersection change shape in the $\zeta$ direction.}
\label{fig:meshintersect}
\end{figure}

\section{Basic properties and consistency of FCIFEM}

Although it is less obvious than in a standard Finite Element formalism, these elements have a partition of unity property, and 
can represent the unity function exactly. Substituting unity in the spline coefficients, rearranging the sums and defining quantities
$\hat{R}_k = \mathcal{R}(R,Z,\zeta,\zeta_k)$ and $\hat{Z}_k = \mathcal{Z}(R,Z,\zeta,\zeta_k)$ (with $\mathcal{Q} = (\mathcal{R},\mathcal{Z},\zeta)$)
yields
\begin{align}
  \sum_k \Omega_{\zeta} \left[ \zeta-\zeta_k \right]
  \sum_{i,k}  \Omega_{R} \left[ \hat{R}_k - R_i \right]
             \Omega_{Z} \left[ \hat{Z}_k - Z_j \right].
\end{align}
Since $\hat{R}_k$ and $\hat{Z}_k$ are independent of $i$ and $j$, the partition of unity property of 1D spline basis functions may be used to
show that the sum over $i$ and $k$ is unity, and finally that the overall expression yields unity.
We do not, however, have the $\delta$ function property that the function value evaluated at a node $\phi(\xb_{i,j,k})$ is equal to the spline coefficient
$\phi_{i,j,k}$ at the node.

The resulting representation is smooth, and will be shown to effectively approximate smooth functions in the large mesh resolution
limit. In the FCIFEM method, even for a simple structured meshing, the domains of support of basis functions associated with two nodes on different surfaces $M$
generally overlap only partially, and in a potentially messy way (see figures \ref{fig:spline_function} and \ref{fig:meshintersect}).
However, given a polynomial mapping function, the representation is piecewise polynomial within a finite set of spatial cells,
with cell faces given by the roots of polynomial equations.

The method is most obviously applicable for anisotropic grids with spacing $\delta \zeta \gg \delta R, \delta Z$, but it is straightforward to map this to
an equivalent problem where the grid spacing is isotropic by compressing the $\zeta$ axis.
In the limit that this isotropic grid is refined equally in each direction, we have a parameter $h = \delta R =\delta Z = \delta \zeta$
representing the grid spacing. 
We wish to show that there exists $\phi \in S$ that is a good approximant to a smooth function
$\bar{\phi}$ with derivatives of order 1 in the $R$ and $Z$ directions, and along the mapping direction
(but which may vary rapidly in the $\zeta$ direction) so that
\begin{equation}
  \left| \phi - \bar{\phi}  \right| < C h^2.
\end{equation}
We will need to assume a certain smoothness of the mapping function so that locally 
$\mathcal{Q} = (\mathcal{R},\mathcal{Z},\zeta_k) = (R_\zeta [\zeta-\zeta_k] + R, Z_\zeta [\zeta - \zeta_k] + Z,\zeta_k) + O(h^2)$ for constants $R_\zeta$ and $Z_\zeta$.
We will also assume that we have at least piecewise linear basis functions so that derivatives exist, and we will take
$\Omega_{R,Z,\zeta} = \Omega$ for simplicity. A constructive proof that good approximants may be found is performed by setting the basis function coefficients
to their nodal values, so $\phi_{i,j,k} = \bar{\phi}(R_i,Z_j,\zeta_k)$. 
In this case we have
\begin{align*}
\phi = \sum_{i,j,k} \bar{\phi}(R_i,Z_j,\zeta_k) & \\
  \times & \Omega \left[\mathcal{R}(R,Z,\zeta,\zeta_k) - R_i \right] \\
  \times & \Omega \left[\mathcal{Z}(R,Z,\zeta,\zeta_k) - Z_j \right]  \\
  \times & \Omega \left[ \zeta,\zeta_k \right] \\
\end{align*}
and we expand about $\xb_0 = (R_0,Z_0,\zeta_0)$, using the near-linearity of the mapping, and
in the region where $\zeta-\zeta_0 = O(h)$, and $|(\mathcal{R}-R_0,\mathcal{Z}-Z_0)| = O(h)$ (which is aligned along the mapping) we find
\begin{align*}  
  \phi   = \sum_{i,j,k} \left(\bar{\phi}(R_0,Z_0,\zeta_0) + h (i,j,k) . \nabla \bar{\phi}(R_0,Z_0,\zeta_0 ) \right) & \\
  \times & \Omega_i \left[R_\zeta (\zeta-\zeta_k) + R \right] \\
  \times & \Omega_j \left[Z_\zeta (\zeta-\zeta_k) + Z \right]  \\
  \times & \Omega_k \left[ \zeta \right] + O(h^2).
\end{align*}
This ordering holds, despite the $\zeta$ derivative of the function being large in general,
because the vector $(i,j,k)$ is aligned almost parallel to the mapping direction for the contibuting basis functions, and the derivatives
are order one in the direction. The algebra proceeds by evaluating the $j$ and $k$ sums, for which $\zeta_k$ is a constant.
From the definition of the derivative of B-Spline, we can show $\sum_i i \Omega_i(x) = x/h $ in the interior region. The splines
also have the partition of unity property $\sum_i \Omega_i(x) = 1 $ in the interior region, so the coefficient of $\bar{\phi}$ is $1$, so 
\begin{align*}
  \phi& = \bar{\phi} - \xb_0.\nabla \bar{\phi} \\
                  & + \frac{\partial\bar{\phi}}{\partial R }     \sum_{k} [R_{\zeta} (\zeta - h k) - R] \Omega_k
                    + \frac{\partial\bar{\phi}}{\partial Z }     \sum_{k} [Z_{\zeta} (\zeta - h k) - Z] \Omega_k
                    + \frac{\partial\bar{\phi}}{\partial \zeta}  \sum_{k} k \Omega_k + O(h^2)\\
    & =  \bar{\phi} + (R-R_0,Z-Z_0,\zeta-Z_0) . \nabla \bar{\phi} + O(h^2)
\end{align*}
and from the smoothness of $\phi$ we have $\phi = \bar{\phi} + O(h^2)$ in the (mapping-aligned) vicinity of any grid point. Similarly,
we can show that the error in derivatives is $O(h)$ in the $R$ and $Z$ directions and along the mapping.

However, better bounds are expected in practice due to an equivalence with more standard discretisations.
For an exact mapping, it is useful to view the discretisation in a `flux tube' type coordinate scheme
$(R',Z',\zeta) = (\mathcal{R}[\xb,\zeta_0],\mathcal{Z}[\xb,\zeta_0],x_\zeta)$. The representation on each surface $M_i$ in these coordinates
with $i\neq 0$ is smoothly distorted by the mapping function, which varies on system scale lengths, so the representational
power for smooth functions is equivalent to the undistorted mapping: we expect it to be consistent to the same order as the original
mapping. The overall representation is then a tensor product of spline functions along the $\zeta$ direction with these $n$th order
consistent representations based on surface $M_i$. At lowest order, the distortion is just a translation in $\mathcal{Q}$ in which case
polynomials of order $n$ are exactly represented on each plane $M_i$. 

\section{Variational formulation} 

Various functional equations for spatial unknowns (typically differential and integro-differential equations) may be represented in
this method via their weak form. For an equation $A(\phi) = 0$, we require that the discrete representation $\phi$, for all weight
functions $\psi$ in the reduced function space, satisfies
\begin{equation}
    \int \psi A(\phi) = 0,
\end{equation}
and leads to a sparse matrix problem $\mathbf{A} \mathbf{k} = 0$ where $k$ represents the coefficients of the FCIFEM representation.

Where $A$ is a local function of $\phi$ (usually a differential operator) the integration involves evaluations of $\phi$ and
$\psi$ at the same spatial locations, and in the standard finite element formalism, the spatial mesh of elements forms a natural
basis for a quadrature (integration) mesh; for polynomial basis elements, appropriate quadratures (such as Gauss points) are
well-known, which allow machine-precision evaluation of these integrals at reasonable cost. In general mesh-free methods, the
lack of alignment between the basis function support domains means that the favourable convergence properties of Gauss quadrature
cannot generally be expected\cite{meshfreequadrature}.

In the set of coordinates $(\mathcal{R},\mathcal{Z},\zeta)$, with an invertible mapping, the domains of integration
are given the tensor product of areas $I$ in the $( \mathcal{R},\mathcal{Z})$ plane (an example of the shapes of such areas is shown in
figure \ref{fig:meshintersect}) multiplied by intervals in $\zeta$. The areas
are bounded by curves which lie on the union of the images of the meshes of neaxby mother planes. For practical examples the boundaries
of these areas are approximately polygonal, with low curvature edges. We do not attempt to do so here, but
calculation of these intersections could be performed with standard packages at least where edges can be taken to be straight.
Boundary conditions would complicate the meshing process substantially. Instead of performing this complicated meshing procedure
we propose the use of simpler methods, in the general spirit of avoiding geometrical complexity.

Increasing the number of Gauss points over that required for a standard mesh problem has been found to be sufficient in tests of certain mesh-free methods\cite{Belytschko1996}. 
We are frequently interested in cases where the operator A is non-local (for example integro-differential, rather than simply differential),
so that evaluation points of $\phi$ and $\psi$ are different\cite{Dominski,Mishchenko} and the usual quadrature approaches are not well-justified.
In the example of this paper we use quadratures specified
on a mesh which conforms to the boundaries\cite{Belytschko1996}; this is only straightforward if the geometry of the boundaries is relatively simple.

\section{Handling boundaries in FCIFEM}

Boundary conditions may be handled in a FEM by generating a mesh that conforms to the boundary surface, and ensuring the finite element
basis functions are consistent with the boundary conditions. However, the mesh volumes in FCIFEM are strongly curved along the $\zeta$ direction, so even
if the boundary is flat in Cartesian coordinates, it is curved in the natural mesh coordinates. Attempting to adapt the FCIFEM meshes so they conform
to these curved surfaces would lead to a somewhat messy meshing problem. For example, the intersection of the mesh volume with the internal region
is in general a complex shape, and in general volumes would need to be broken into smaller pieces to simplify the geometry.
In conjunction with this, an additional coordinate transform would in general need to be used to map the curved surfaces to flat faces via an
isoparametric transform. The philosophy here, and for mesh-free meshods in general, is to avoid these geometrical complications.

There are a number of approaches to handling boundaries in mesh-free methods\cite{Chavanis2005,Huerta2004}. In the methods close in spirit to that proposed here,
essential boundary conditions can be imposed by transforming shape functions so that they conform, by introducing penalty functions to 
the minimisation problem resulting from the weak form, or by introducing an additional boundary mesh that conforms exactly. We will use the latter
method, where a standard finite element method mesh is defined near the boundary.

If we denote a function represented by the FCIFEM representation in eq \ref{eq:FCI_FEM} as $\phi_{FCI}$ and a conventional finite element representation as
$\phi_{FEM}$, with
\begin{equation}
  \phi_{FEM} = \sum_{i} c_i K_i(\xb),
\end{equation}
we can produce a blended representation
\begin{equation}
 \phi_{BLE} = B(\xb) \phi_{FEM}(\xb) + [1-B(\xb)] \phi_{FCI}
\end{equation}
where $B(\xb)$ is a {\it ramp function} with $B=1$ on the boundary $\delta \Sigma$, and $B=0$ in the bulk of the domain $\Sigma$ apart from a
narrow region near the boundary, where the function smoothly ramps from $1$ to $0$. The FEM representation is taken to conform to the boundary,
so the boundary condition is exactly satisfied in the appropriate restriction of the FEM function space. 

Also, since the FCIFEM and FEM representations individually are consistent up to some order, the blended sum of the two representations is also able
to exactly represent polynomials up to that order. However, if the nodal values of the FEM and FCI representation are specified independently,
we have introduced additional degrees of freedom in the overall representation, and these may not be sufficiently orthogonal  
(there may be two ways of approximately representing the same spatial function with very different
sets of coefficients): this may result in an ill-conditioned or singular matrix problem. This difficulty can be dealt with in general by modifying
the shape functions of the representation in the boundary region to ensure orthogonality\cite{Belytschko1995,Huerta2004}, but for the more specialised representation
chosen here, there is a simpler way to proceed. We require that where the FEM representation and the FCIFEM representation share nodal positions, 
they have the same nodal values $c_i = \phi_j$; this reduces the number of extra
degrees of freedom, but still ensures consistency.

For the regular FCIFEM grids used in this paper's examples, the FEM boundary grid will be specified on the same $(R,Z)$ node positions; because the FEM
grid does not directly capture the anisotropy, the spacing in the $\zeta$ direction needs to be finer than the FCIFEM grid to fully capture the $\zeta$
depencence near the boundary. We have, however, used the same $\zeta$ grid for the FEM and FCIEFEM in the example problems below. The ramp function is
chosen to be the sum of FEM basis functions associated with boundary nodes\cite{Belytschko1995}.


\section{A simple 2D periodic example problem}

In order to provide a simple test case, as well as a straightforward comparison against earlier methods for representing
anisotropic functions, we consider a problem in a doubly periodic domain $Z \in [0, 2 \pi]$ and $\zeta \in [0, 2 \pi]$ with a mapping induced by a
straight field $\Bb = \zetahat B_{\zeta} + \Zhat B_Z $ and $B_{\zeta} = B_z = 1$. The differential equation we consider in the
remainder of the paper (typical for the problems of electromagnetic coupling which arise in tokamak turbulence) is the Laplacian inverse
\begin{equation}
   \nabla^2 \phi = \rho
\end{equation}
where we solve for $\phi$. We choose
\begin{equation}
  \rho(Z,\zeta) = \sin[ n ( Z - \zeta) ] [1+\sin( Z )]/2,
\end{equation}
with $n=10$, which results in a nearly field aligned perturbation, since the spatial wavenumbers of the perturbation
are $(k_Z,k_{\zeta}) = (n,-n)$ and $(n \pm 1,-n)$, which are aligned with the field to within $10\%$.
Note that the periodicity of the grid and homogeneity of the
problem would allow direct use of discrete Fourier space to solve for the spline coefficients; the analytic solution is
also easily obtained using the Fourier method.

For this perturbation, which has an anisotropy direction aligned $pi/4$ radians from the $Z$ and $\zeta$ direction, a Cartesian
tensor spline representation on the $Z$ and $\zeta$ directions needs to have high resolution in both directions, and it is most
efficient to choose similar resolution in $Z$ and $\zeta$ directions. On the other hand, the FCIFEM representation can take
advantage of the slower spatial variation along the anisotropy.

\begin{figure}[htb]
\centering
\includegraphics[width=8cm]{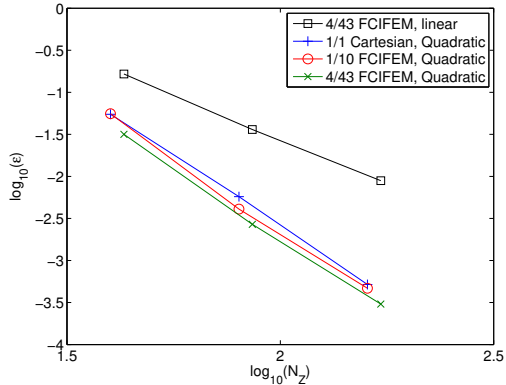}
\caption{  Convergence of the $L^2$ error for a simple analytical test problem in a 2D periodic domain with the FCIFEM method
  and a straight-line mapping function. The legend gives the ratio of $N_{\zeta}/N_Z$ grid points for each scan. All
  these testcases use quadratic splines except the uppermost data. The data marked with crosses uses a standard cartesian
  spline representation rather than FCIFEM.
   }
\label{fig:conv_simple}
\end{figure}

For this simple case, if $N_Z$ is an integer multiple of $N_{\zeta}$, the FCIFEM representation reduces to a standard tensor
product spline representation, in a sheared coordinate system, because the mapping between surfaces $M$ aligns the nodes.
A scan is performed over resolution for both $N_{\zeta}/N_Z = 4/43$ and $N_{\zeta}/N_Z = 1/10$ in order to demonstrate that the
accuracy of the method is not significantly degraded by lack of node alignment. The error a a function of grid resolution is
plotted in figure \ref{fig:conv_simple}. The testcases for $N_{\zeta}/N_Z = 4/43$ are repeated for linear splines,
but other results are reported for quadratic splines only.
We also compare the results with the spline representation on a regular Cartesian grid in $Z$ and $\zeta$
(with $N_{\zeta}/N_Z=1$) to demonstrate the advantage of the anisotropy-capturing scheme; the results of the
anisotropic scheme have similar error for the same resolution in $N_Z$, but require a factor of 10 fewer
grid points in total. RMS errors converge as very nearly $1/N^2$ for linear splines, and $1/N^3$ for quadratic splines.

\section{A tokamak-related example problem}

\begin{figure}[htb]
\centering
\includegraphics[width=8cm]{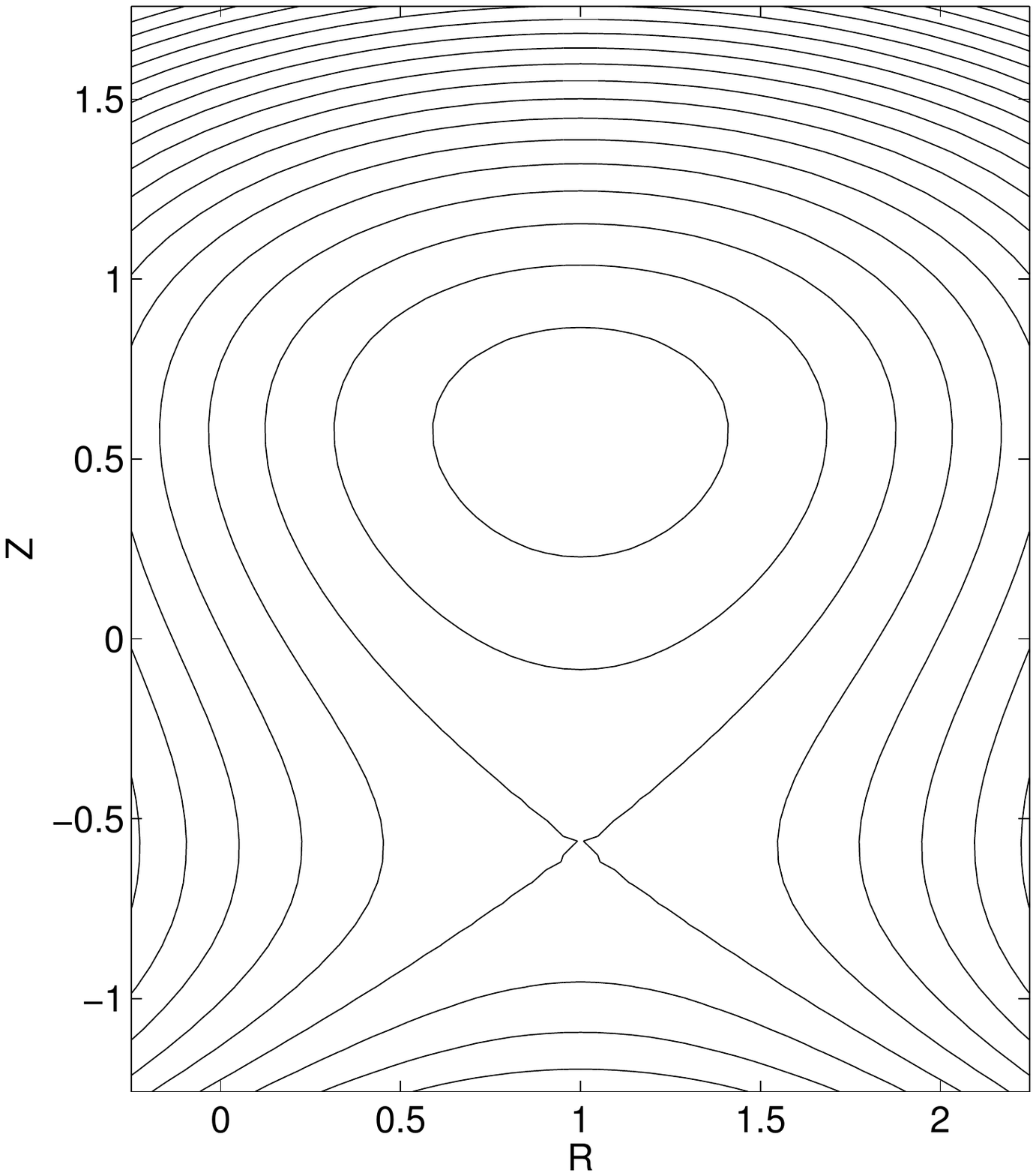}
\caption{ The  $\Bb$ field lines on the
  $(R,Z)$ plane (the $\zeta$ direction is into the page); these are also
  the countours of the function $A$ specified in equation \ref{eq:aeq}.
   }
\label{fig:figacont}
\end{figure}

The motivating problem is magnetic confinement fusion; we consider a plasma is confined by a magnetic field, with the field lines shown
in fig. \ref{fig:figacont}. If the bounding rectangle is taken to be a
physical wall, the volume can be separated into an `open field line' region, whose magnetic field lines intersect the wall,
and a closed field line region. Fig. \ref{fig:figacont} is typical of a `diverted' tokamak configuration, where the nested
set of flux surfaces end at a separatrix, and points outside the last closed flux surface are connected by the magnetic
field lines to the wall.

A common computational task is to trace the trajectories of particles in a turbulent field generated by some set of charges and currents,
which may be part of a particle-in-cell simulation\cite{BirdsallLangdon,SebORB5}. Accuracy of particle tracing requires a smooth
representation of the fields, and a partial differential equation (in general integro-differential) must be solved to determine the
fields based on currents and charge sources; we will explain how to use the partially mesh free method to solve an example problem
of this type. 

We define a magnetic field
\begin{equation}
   \mathbf{B} = B_0 \boldsymbol{\hat{\zeta}}  + \boldsymbol{\hat{\zeta}} \times \nabla A(R,Z)
\end{equation}
with
\begin{equation}
  A(R,Z) = (R-1)^2 + Z ( Z^2 - 1 )
  \label{eq:aeq}
\end{equation}
which represents a diverted configuration in the large aspect ratio limit. The contours of $A$ (figure \ref{fig:figacont})
are the field lines projected onto the $(R,Z)$ plane, and an X-point is seen at $(R,Z) = (1, -1 / \sqrt{3} )$.

The mapping function is approximated by using a Taylor series expansion of the mapping function $\mathcal{Q}(\xb,\zeta)$ up to
quadratic order in $\zeta$ and a spline representation of the Taylor series coefficients on the $(R,Z)$ grid used for the FCIFEM
elements. We quantify the error by tracing accurate and approximate field lines starting at positions in the $(R,Z)$ spatial domain
at $\zeta=0$ to $\zeta=\zeta_1$. For field lines which remain in the simulation volume, the RMS error in the final $(R,Z)$ position
is $0.016$, and the maximum error is $0.06$ (the error is concentrated near the boundaries at large $Z$ where $A$ has large gradients). 


As in the previous section, we solve the Laplacian inverse problem, but now over a rectangular domain in $R,Z$ and a periodic $\zeta$ direction
with $\phi=0$ on the boundaries of the domain. We take $R \in [R_0,R_1]$ and $Z \in [Z_0, Z_1]$ and $\zeta \in [0,\zeta_1]$.
In order to ensure that the representation satisfies the boundary conditions exactly, we use a single cell-width first order
standard FEM formalism near the boundary, blended with the FCIFEM method using a ramp function. The ramp function is taken to be
the sum of the first-order FEM basis functions, which is simply a linear function $(1-x)$, where $x$ is the distance to the
boundary in grid units, for points away from a corner; for points near a corner the ramp function is $(1-x)(1-y)$ with $x$ and $y$ the distance
to the nearby boundary edges in grid units.


Quadratic B-Splines basis functions are chosen to represent the field, with uniform grid spacing $\delta R$, $\delta Z$, $\delta \zeta$.

Discretisation proceeds by taking the weak form and finding $\phi \in S $ such that
\begin{equation}
  \int dV ( \nabla \psi . \nabla \phi + \psi \rho) = 0 
\end{equation}
for all $\psi \in S$. The integration is performed using a set of quadrature points evenly spaced in $R,Z$ and $\zeta$, 10 times finer
than $\delta R$, $\delta Z$, $\delta \zeta$ respectively.

\subsection{A basic convergence test}

In order to examine the basic convergence of the method, a simple test problem is considered with $\rho(R) = \sin(0.5 \pi R) \sin(\pi [Z+1.0]/2.5)$.
the domain $R_0 = 0$,  $R_1 = 2$,  $Z_0 = -1$,  $Z_1 = 1.5$ and $\zeta_1 = \pi/20$ is chosen. The solution $\phi$ to this test-problem is not
aligned along the field, but constant along $\zeta$, so the use of the FCIFEM is not advantageous in this case.
Along the field line, the perturbation varies with typical scale length $|B_{\zeta}|/|B_{R,Z}|$
times longer than typical wavelengths in $R$ and $Z$. We have therefore chosen $\delta Z$ comparable to $ \delta R |B_{\zeta}|/|B_{R,Z}|$
so that the effective resolution is sufficiently high along the field line.

\begin{figure}[htb]
\centering
\includegraphics[width=8cm]{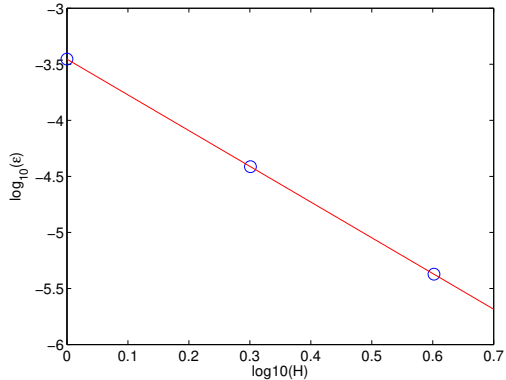}
\caption{ Convergence of the $L^2$ error for a simple analytical test problem in a 3D domain rectangular in $R,Z$ and periodic in $\zeta$
          with the FCIFEM method and a mapping function based on a divertor-type field. }
\label{fig:convergence}
\end{figure}

The field-aligned mesh leads to projected domains of the basis functions in the $(R,Z)$ plane
of extent $\sim |\mathcal{Q}[R,Z,\delta \zeta]-\mathcal{Q}[R,Z,0]| + |\delta R|+|\delta Z|$, so the representation, and the
effective resolution in the $R$ and $Z$ direction depend on the number of $\zeta$ points chosen. The number of grid points
at lowest resolution $(N_R,N_Z,N_{\zeta}) = (20,20,1)$, and this is uniformly increased by a factor $H$ in each direction to perform
a convergence scan. Figure
\ref{fig:convergence} shows the $L^2$ error $\epsilon$ in the solution versus $H$, which drops as $ H^{3.2} $, in line with
the value $3$ expected for a standard second order FEM, and as good as could be expected for the degree of representation smoothness
chosen.

\subsection{An illustrative tokamak problem}

In order to illustrate the treatment of anisotropic structures in this method, we consider a second test problem 
in the same spatial domain as the previous problem.
Here, we take $\rho(\Rb) = \int d\zeta [ \delta(\mathbf{R} - \xb(\zeta)) - \delta(\mathbf{R} - \xb(\zeta) - \hat{\zeta} \zeta_1/2) ]$
with $\xb(\zeta)$ the curve traced out by a field
line parameterised by $\zeta$, starting from the point $\xb(0) =(0.36,-1.0,0)$, so that it passes near the X-point. This results in a
highly anisotropic charge perturbation along the magnetic field lines with neighbouring charge filaments of opposite sign, typical of a localised unstable
drift mode. The number of grid points in each direction is $(N_R,N_Z,N_{\zeta}) = (100,100,1)$; because the perturbation is highly field aligned, the solution
is expected to be also well-aligned, so that the anisotropic solution is well captured, despite using only $1$ point in the $\zeta$ direction.
The mean anisotropy can be quantified by the RMS angle that the field line makes with respect to the $\zeta$ direction: following a field line from $\zeta_0$ to $\zeta_1$,
the RMS displacement on the $(R,Z)$ plane is $0.3$ units, corresponding to 16 grid cells. Thus, to capture the anisotropy using a Cartesian grid, we would require roughly
$16$ grid cells in the $\zeta$ direction.

The resulting matrix problem, coupling the coefficients of the $\phi$ representation to those of the weight function, is a sparse matrix
of rank equal to the number of degrees of freedom of the system, which we treat as being unstructured. Each row in the matrix has of order
$200$ non-zero entries when quadratic splines are used; this is somewhat less sparse, due to irregular overlapping of domains of support,
than in the corresponding 3D tensor spline representation, where $125$ non-zero entries would be expected.
Index reordering is quite effective in
reducing the bandwidth of the resulting matrix, so that direct solution using banded matrix calculation is straightforward for the test
problem under consideration. 

\begin{figure}[htb]
\centering
\begin{tabular}{@{}l@{}}
\subfigure[]{
\includegraphics[width=8cm]{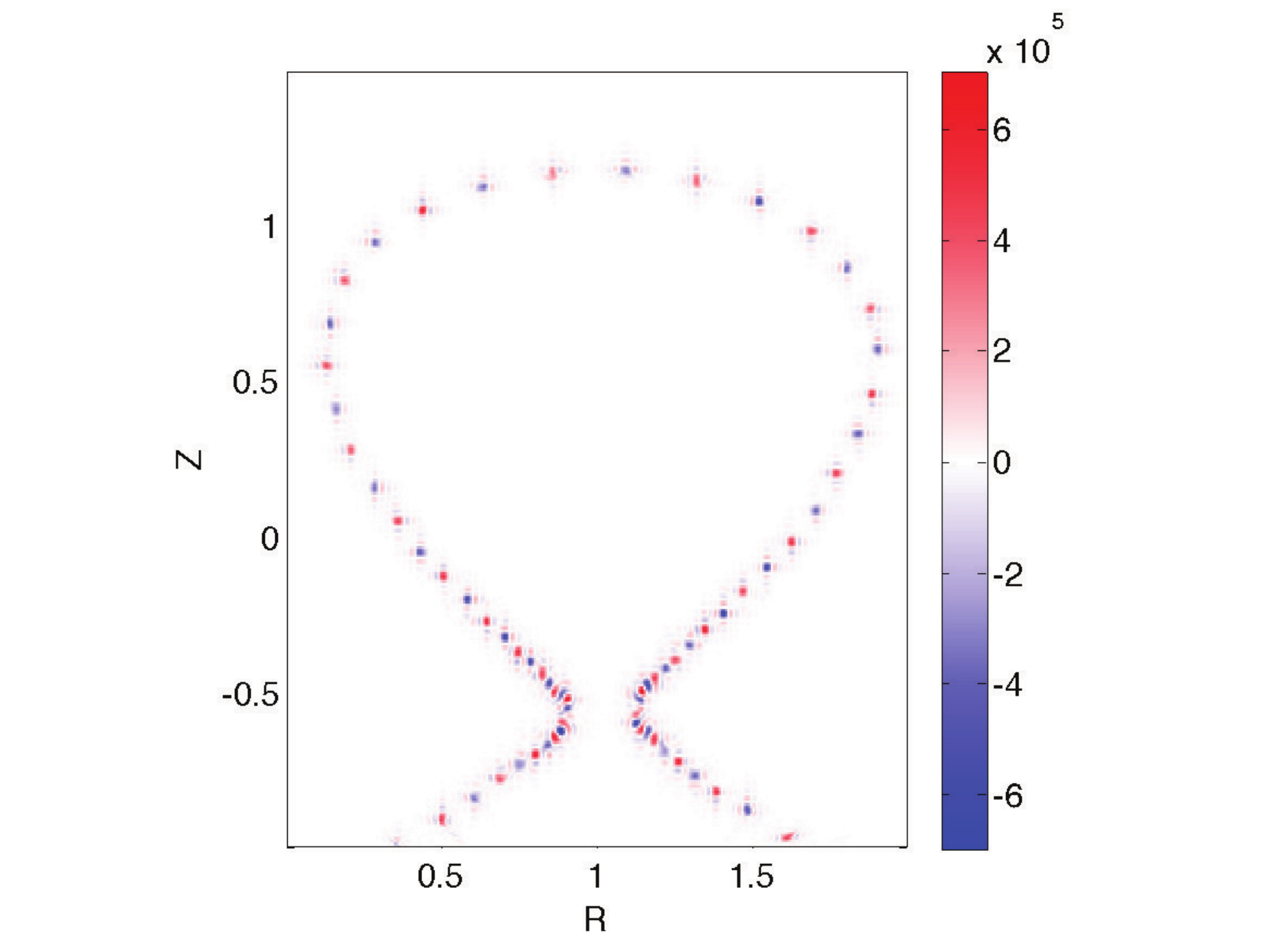}
\label{fig:rhobar_2d}
}
\subfigure[]{
\includegraphics[width=8cm]{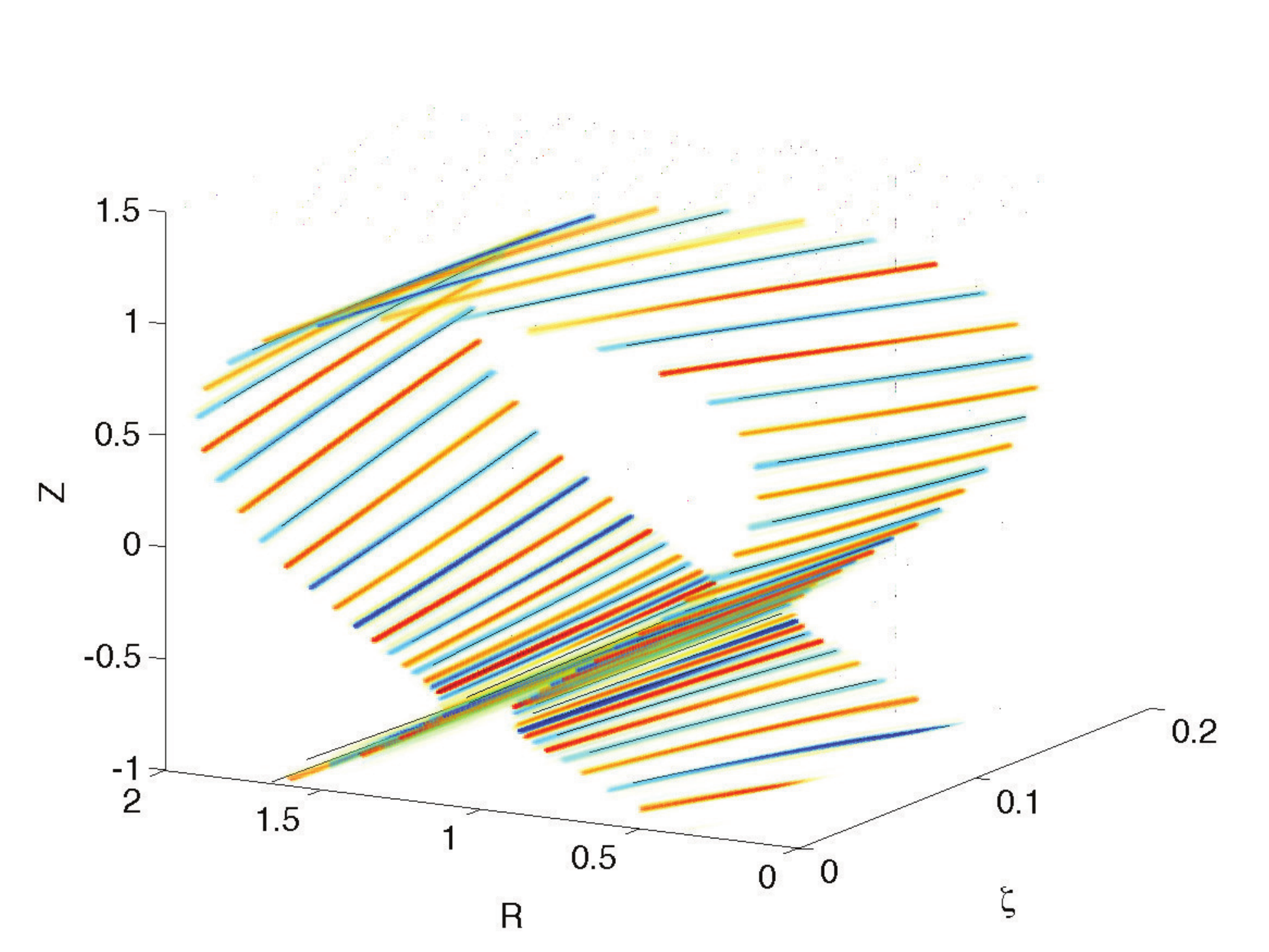}
\label{fig:rhobar_3d}
} \\
\subfigure[]{
\includegraphics[width=8cm]{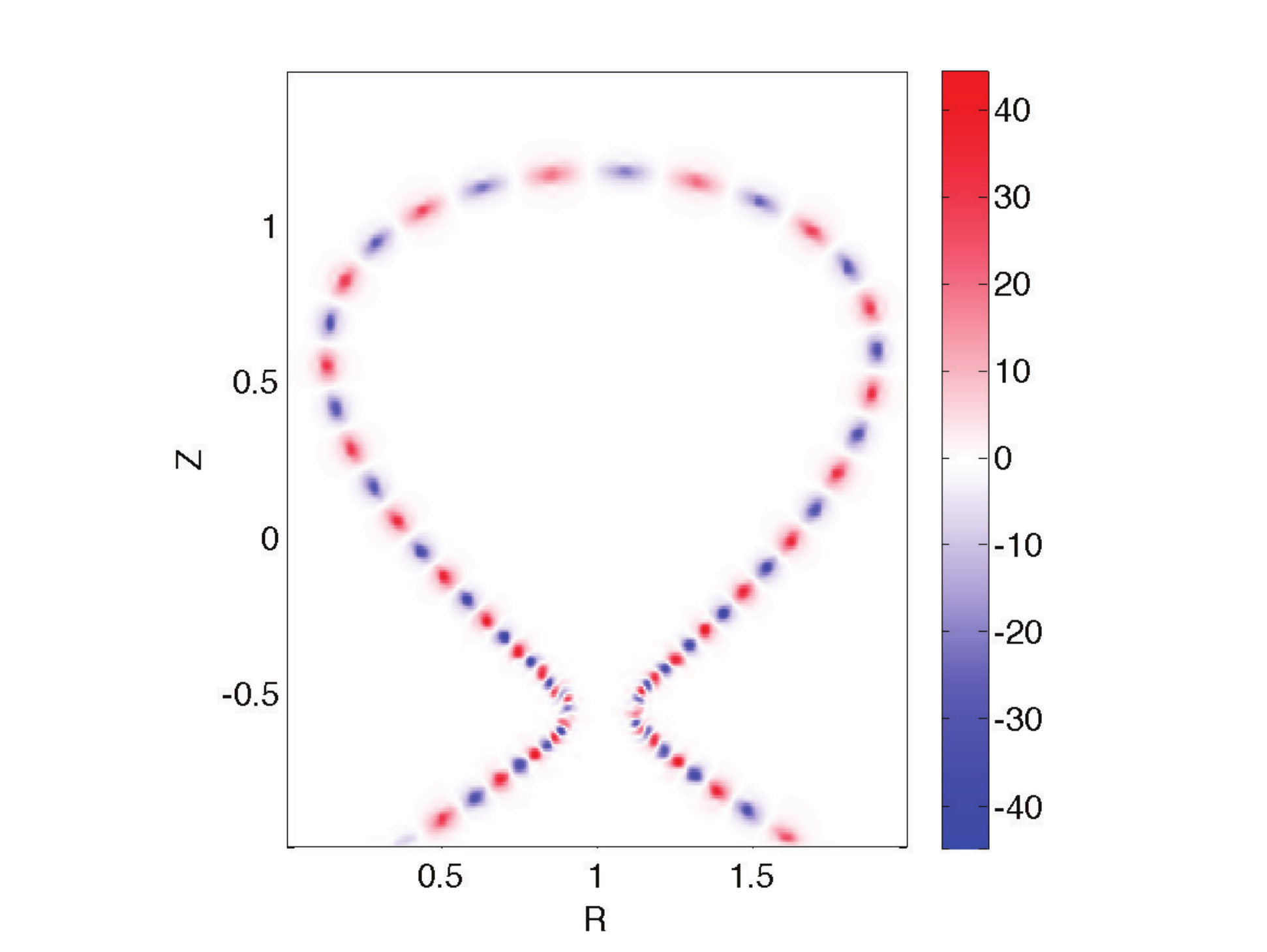}
\label{fig:phisoln_2d}
}
\subfigure[]{
\includegraphics[width=8cm]{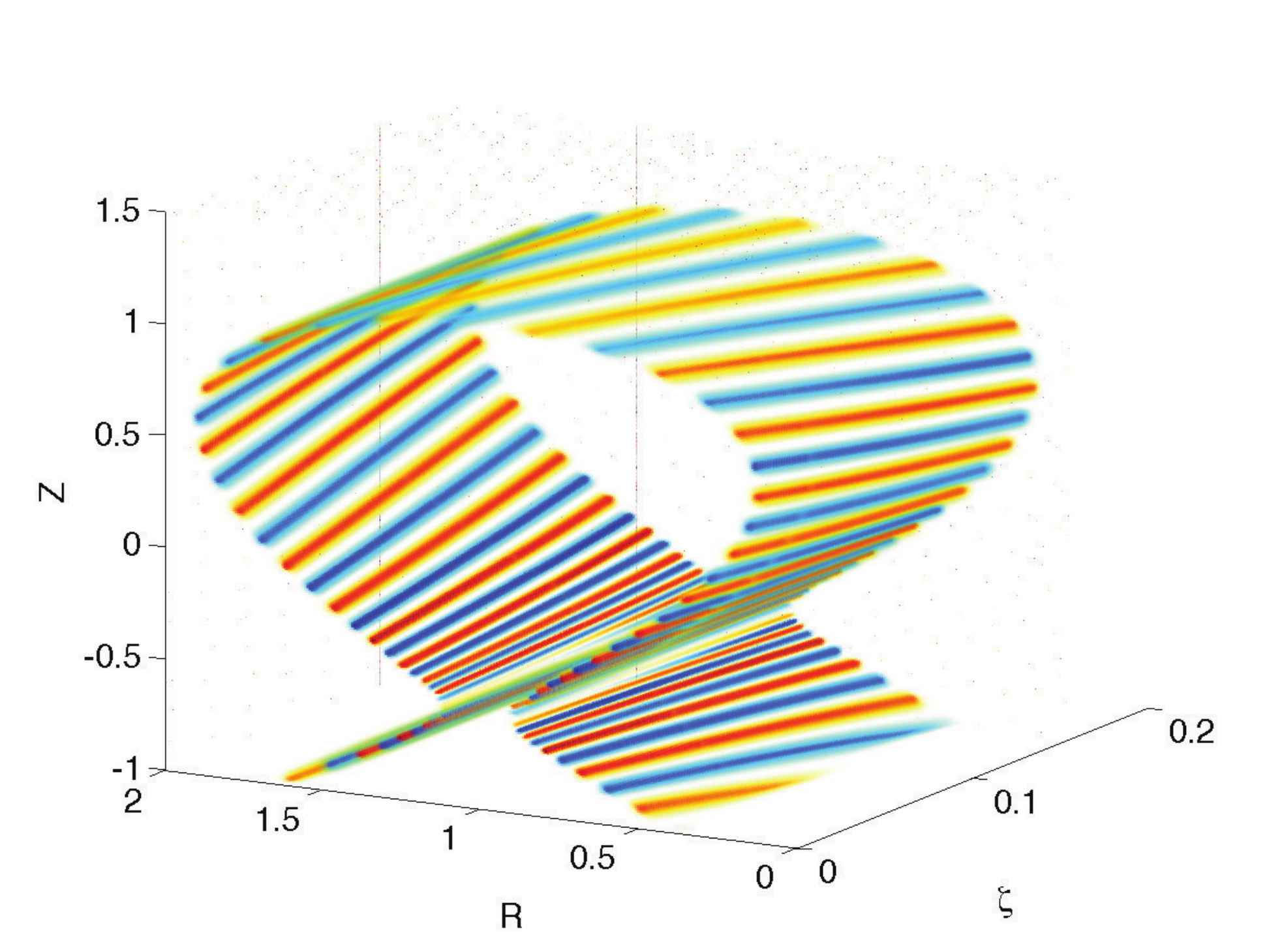}
\label{fig:phisoln_3d}
}
\end{tabular}
\caption{ 2D plots versus $R$ and $Z$ at $\zeta=0$ (a,c), and 3d voxel renderings (b,d), of $\bar{\rho}$ (a,b) and $\phi$ (c,d). Plot
  (b) shows a field line in black, essentially coincident with the negative charge region. 2D plots are oversampled
  by a factor of three compared to the spline grid. 3D plot show regions of positive charge/potential as red, and negative as blue.
  The $\zeta$ direction extends periodically with period $\pi/20$, but
  repeats are not shown; note that the scales are not equal in the 3D plots. }
\label{fig:3dfield}
\end{figure}

To show the projection of the charge $\rho$ into the space $S$, the weak form
\begin{equation}
  \int dV ( \psi \bar{\rho} - \psi \rho) = 0 
\end{equation}
is solved for $\bar{\rho} \in S$ in the same fashion as for the Laplacian problem. We show 2D plots and 3D plots of $\bar{\rho}$ and
$\phi$ in figure \ref{fig:3dfield}. Due to the structure of the charge (alternating charge lines of opposite sign),
the potential $\phi$ decays rapidly away from the field lines where $\rho$ is nonzero. The strong anisotropy of the charge and of
$\phi$ are clear in the 3D plots. The structures in the reproduced field $\bar{\rho}$ are conicident with the field lines, as shown in plot (b);
this provides some evidence that the mapping $\mathcal{Q}$ is sufficiently accurate. Note that the spatial variation of the the imposed $\rho$
becomes too rapid to be captured by the spline representation near the X-point; this is not typical of the turbulent structures
which tend to be of a more uniform typical wavenumber, but allows us to demonstrate how the
scheme behaves at short spatial scale. The integration inherent in this Galerkin type method averages out the positive and negative variations below the grid scale
to zero in this region. The FCIFEM method, which has uniform resolution near the X-point, is better able to resolve structures in this region than
conventional schemes for MCF problems: conventional methods have low resolution near the X-point, as they place sets of nodes on flux surfaces,
and flux surfaces become widely spaced near the X-point.

Behaviour near the boundaries is acceptable from visual inspection in the reproduced field and the solution to the Laplacian problem. In the interior region,
there is some oscillation evident in the function $\bar{\rho}$ in figure
\ref{fig:rhobar_2d} which is attempting to represent a $\delta$ function using a smooth spline representation:
this is standard for finite-element type representations and not a particularity of FCIFEM. As with the behaviour near the X-point,
we have deliberately chosen a somewhat difficult example that probes how the scheme handles sub-grid scale forcing.
To determine the effect of an inexact mapping on the problem, the Laplacian problem was repeated with a precise mapping used based on numerical
integration rather than a low-order Taylor series. Visually the results were indistinguishable. The relative RMS difference in solutions
using the exact and inexact mapping was $3\%$. 

\begin{figure}[htb]
\centering
\includegraphics[width=8cm]{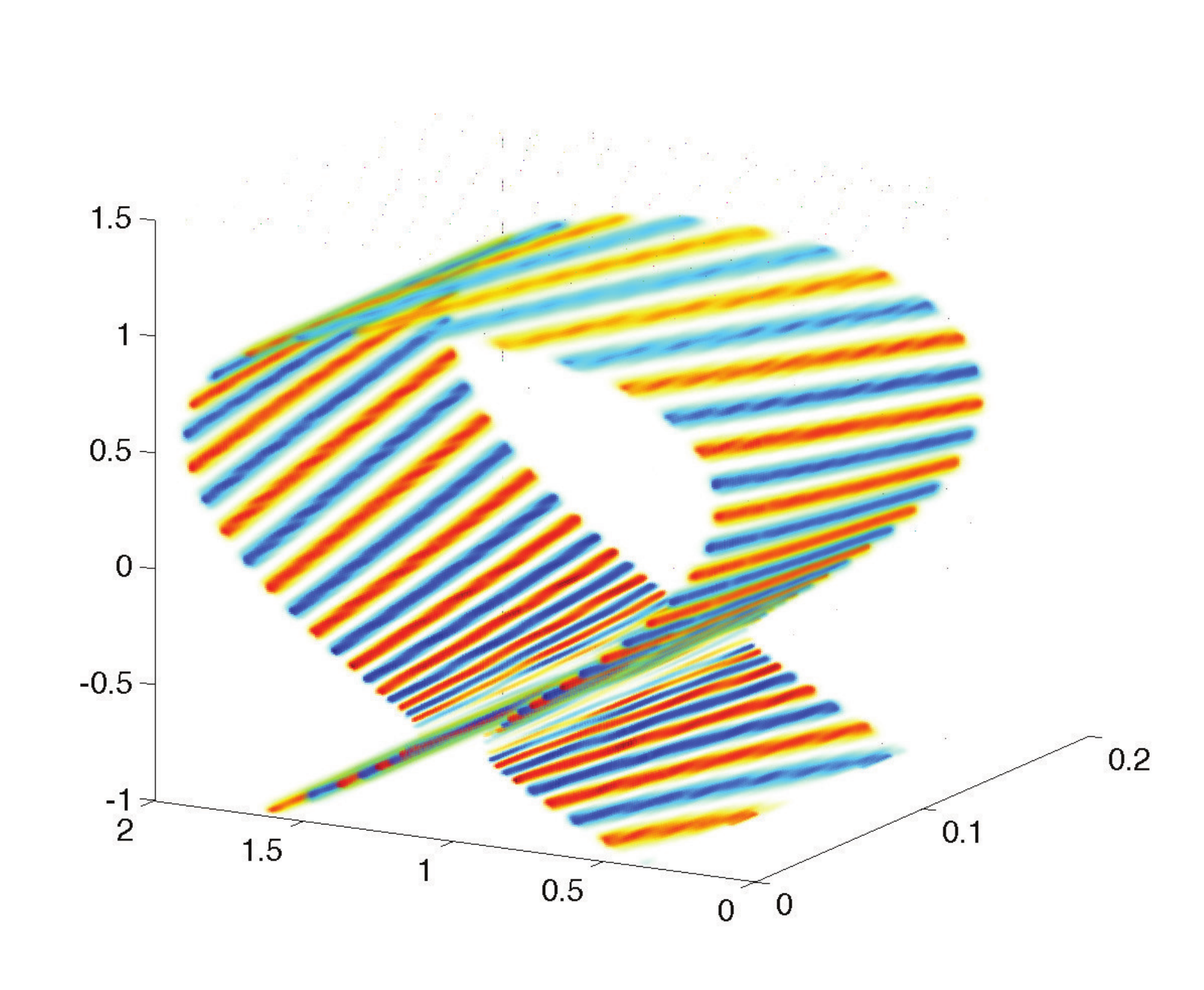}
\caption{3d voxel rendering of $\phi$ for a finite element method based on Cartesian grid.}
\label{fig:3dfield_cart}
\end{figure}

To demonstrate the effectiveness of the FCIFEM method, we compare it to a method based on a Cartesian mesh with $10$ equally spaced
points in the $\zeta$ direction and $100$ points in both $R$ and $Z$ directions; the Cartesian representation has 10 times as many node
coefficients as the FCIFEM representation, and solution of the matrix problem requires an iterative method. The solution $\phi$ plotted in fig. \ref{fig:3dfield_cart} is
visually quite similar to that of the FCIFEM, although it is noticably less smooth along the field lines; there are short wavelength
oscillations at the grid scale which arise as the Cartesian method attempts to reproduce short wavelength anisotropic
structures. This kind of aliasing error may be problematic even if the usual measures of error are small. For example,
where the fields are used to evaluate particle orbits, the short wavelength structures translate into small timesteps;
this may be worked around by filtering out grid-scale wavelengths\cite{SebORB5} at the cost of higher grid resolution.

\section{Discussion}

We have introduced a numerical technique for representing anisotropic functions, and solving equations related to these functions, in regions where the direction of anisotropy
is spatially varying. The guiding principle (as with the FCI method) is to incorporate the complexity of the field line geometry into a mapping function, and avoid
specialised geometrical representations based on the magnetic field topology.

One important feature of the method proposed here is the ability of the method to handle curved anisotropic structures; once the basis and mapping
is defined, there is a straightforward and systematic method for evaluating differential operators. This has been handled only in part in earlier
methods. The difficulty is that in MCF problems, the anisotropic stuctures extend along field lines to a length scale typically of order the system scale $R_g$,
and the departure of field lines from straight lines over this scale is also of this order, even in cylindrical coordinates
(the departure from straightness is of order $a$, the {\it minor radius}, in a tokamak problem). The wavelengths of turbulent structures perpendicular to the
field, on the other hand, are orders of magnitude smaller. Methods that require the anisotropy direction to be constant in each mesh cell will require a finer
mesh spacing along the field line that those able to explicitly incorporate curvature, like the FCIFEM.

We have chosen to test an inexact mapping for the tokamak geometry testcase in the previous section. For this special case, it is possible to take advantage of the
topology of the problem to produce a near-exact numerical approximation to the field line mapping operation, but the philosophy of this approach is to avoid relying on
concepts like flux coordinates that would be used to construct such a map. Since the method appears to be quite robust to the use of even quite crude inexact mappings,
we expect the choice to be a matter of convenience. 

\section*{Acknowledgments}
This work has been carried out within the framework of the EUROfusion Consortium and has received funding from the Euratom research and training
programme 2014-2018 under grant agreement No 633053. The views and opinions expressed herein do not necessarily reflect those of the European Commission.
Thanks to Eric Sonnendrucker for helpful comments on a draft version of this document.

\bibliography{gyrokin}

\end{document}